\documentclass[5pt]{article}
\usepackage{amsmath}
\usepackage[latin1]{inputenc}
\usepackage{amsmath,amsthm,amssymb}
\usepackage{graphicx}
\usepackage{color}

\renewcommand{\theequation}{\thesection.\arabic{equation}}

\newtheorem{lemma}{Lemma}[section]
\newtheorem{thm}{Theorem} [section]

\newtheorem{rem}{Remark}[section]
\title{ Diophantine equations involving Euler function
\thanks {The work is supported by the Natural Science Foundation of China with No.11671153, No.11971180.}}
\author{{Hairong Bai\footnote{\ Hairong Bai, Institute of Mathematics Science, South China Normal University,
 \small Guangzhou, Guangdong, China (Email:baihairong2007@163.com)}
}\ }
\date{}
\begin{document}
\baselineskip15pt \maketitle
\renewcommand{\theequation}{\arabic{section}.\arabic{equation}}
\catcode`@=11 \@addtoreset{equation}{section} \catcode`@=12
\begin{abstract}

In this paper, we show that the equation $\varphi(|x^{m}-y^{m}|)=|x^{n}-y^{n}|$ has no nontrivial solutions in integers $x,y,m,n$ with $xy\neq0, m>0, n>0$ except for the solutions $(x,y,m,n)=((2^{t-1}\pm1),-(2^{t-1}\mp1),2,1), (-(2^{t-1}\pm1),(2^{t-1}\mp1),2,1),$ where $t$ is a integer with $t\geq 2.$ The equation $\varphi(|\frac{x^{m}-y^{m}}{x-y}|)=|\frac{x^{n}-y^{n}}{x-y}|$ has no nontrivial solutions in integers $x,y,m,n$ with $xy\neq0, m>0, n>0$ except for the solutions $(x,y,m,n)=(a\pm1, -a, 1, 2), (a\pm i, -a, 2, 1),$ where $a$ is a integer with $i=1,2.$


 \end{abstract}

 \section{Introduction}
 $ $

Euler function is defined as $\varphi(n)=\sharp\{r: r\in Z, r>0, (r, n)=1\}.$ Many famous problems on Euler function $\varphi$ have been studied.
As it is well known, there are many Diophantine equations involving Euler function. For example the equation $\varphi(m)=\varphi(n)$ (see [2], [18], [19], [32]), $k\varphi(n)=n-1, \varphi(m)=\sigma(n),$  the iteration problem of functions $\varphi$ and $\sigma,$ and so on.

In 2005, Luca[24]proved that, if $b\geq 2,$ is a fixed integer, then the equation
$$\varphi(x\frac{b^{m}-1}{b-1})=y\frac{b^{n}-1}{b-1}, x,y \in\{1,2,\cdots, b-1\}$$
has only finitely many positive integer solutions$(x,y,m,n).$

In 2017,Yong-Gao chen and Hao tian [8] proved that,
The equation $$\varphi(x^{m}-y^{m})=x^{n}-y^{n}$$ has no solutions in positive integers $(x,y,m,n)$ except for the trivial solutions $(x,y,m,n)=(a+1,a,1,1),$ where $a$ is a positive integer.
The equation $$\varphi(\frac{x^{m}-y^{m}}{x-y})=\frac{x^{n}-y^{n}}{x-y}$$ has no solutions in positive integers $(x,y,m,n)$ except for the trivial solutions $(x,y,m,n)=(a,b,1,1)$ where $a,b$ are integers with $a>b\geq 1$

In this paper, Using the similar method of [8], we consider the following equations.
$$\varphi(|x^{m}-y^{m}|)=|x^{n}-y^{n}|, \eqno{(1.1)}$$
$$\varphi(|\frac{x^{m}-y^{m}}{x-y}|)=|\frac{x^{n}-y^{n}}{x-y}|, \eqno{(1.2)}$$
in integers $x, y, m, n$ with $xy\neq 0, m>0, n>0.$

The equation (1.1) has the trivial solution in integers $(x,y,m,n)=(a\pm1, a, 1, 1),$ where $a$ is a integer. The equation (1.2) has the trivial solution in integers $(x,y,m,n)=(a, b, 1, 1), (\pm1,\mp1,2r+1,2\mu+1),$ where $r, \mu, a, b$ are integers.

In this paper the following results are proved.

\begin{thm}\label{thm1}
The equation (1.1) has no nontrivial solutions in integers $(x,y,m,n)$ with $xy\neq 0, m>0,n>0$ except for the solutions $(x,y,m,n)=((2^{t-1}\pm1),-(2^{t-1}\mp1),2,1), (-(2^{t-1}\pm1),(2^{t-1}\mp1),2,1)$ where $t\geq 2$ is a integer with $t\geq 2.$
\end{thm}

\begin{thm}\label{thm1}
The equation (1.2) has no nontrivial solutions in integers $(x,y,m,n)$ with $xy\neq 0, m>0, n>0$ except for the solutions
$(x,y,m,n)=(a\pm1, -a, 1, 2), (a\pm i, -a, 2, 1),$ where $a$ is a integer with $i=1,2.$
\end{thm}

we always assume that $|x|>|y|\geq 1, m,n\geq 0, m\neq n$ without loss of generality. According to the positive, negative of $x, y$ and the parity of $m, n,$, the equation (1.1) is equivalent to the following equations(a1-a4):
$$\varphi(|x|^{m}-|y|^{m})=|x|^{n}-|y|^{n}, xy<0, 2|m, 2|n \quad or \quad xy>0 \eqno{(a1)}$$
$$\varphi(|x|^{m}+|y|^{m})=|x|^{n}+|y|^{n}, xy<0, 2\nmid mn \eqno{(a2)}$$
$$\varphi(|x|^{m}-|y|^{m})=|x|^{n}+|y|^{n}, xy<0, 2|m, 2\nmid n \eqno{(a3)}$$
$$\varphi(|x|^{m}+|y|^{m})=|x|^{n}-|y|^{n}, xy<0, 2\nmid m, 2|n \eqno{(a4)}$$

The equation (1.2) is equivalent to the following equations(a5-a9):
$$\varphi(\frac{|x|^{m}-|y|^{m}}{|x|-|y|})=\frac{|x|^{n}-|y|^{n}}{|x|-|y|}, xy>0 \eqno{(a5)}$$
$$\varphi(\frac{|x|^{m}-|y|^{m}}{|x|+|y|})=\frac{|x|^{n}-|y|^{n}}{|x|+|y|}, xy<0, 2|m, 2|n \eqno{(a6)}$$
$$\varphi(\frac{|x|^{m}+|y|^{m}}{|x|+|y|})=\frac{|x|^{n}+|y|^{n}}{|x|+|y|}, xy<0, 2\nmid mn \eqno{(a7)}$$
$$\varphi(\frac{|x|^{m}-y^{m}}{|x|+|y|})=\frac{|x|^{n}+|y|^{n}}{|x|+|y|}, xy<0, 2|m, 2\nmid n \eqno{(a8)}$$
$$\varphi(\frac{|x|^{m}+y^{m}}{|x|+|y|})=\frac{|x|^{n}-|y|^{n}}{|x|+|y|}, xy<0, 2\nmid m, 2|n \eqno{(a9)}$$

In summary, we just need to consider the following equations (1.3-1.6) in positive integers $x,y,z,m,n$ with $x>y\geq1, m\neq n.$
$$\varphi(z\frac{x^{m}-y^{m}}{x-y})=z\frac{x^{n}-y^{n}}{x-y} \eqno{(1.3)}$$
$$\varphi(z\frac{x^{m}+y^{m}}{x+y})=z\frac{x^{n}+y^{n}}{x+y}, 2\nmid mn \eqno{(1.4)}$$
$$\varphi(z\frac{x^{m}-y^{m}}{x+y})=z\frac{x^{n}+y^{n}}{x+y}, 2\mid m, 2\nmid n \eqno{(1.5)}$$
$$\varphi(z\frac{x^{m}+y^{m}}{x+y})=z\frac{x^{n}-y^{n}}{x+y}, 2\nmid m, 2|n \eqno{(1.6)}$$

In fact equations (a1-a4) and (a5, a7-a9) are special cases for $z=x+y$ and $z=1,$ respectively. The equation (a6) is equivalent to the equation $$\varphi((|x|-|y|)\frac{|x|^{m}-|y|^{m}}{|x|^{2}-|y|^{2}})=(|x|-|y|)\frac{|x|^{n}-|y|^{n}}{|x|^{2}-|y|^{2}}, xy<0, 2|m, 2|n.$$
By [9], the equation (a1, a5, a6) has no nontrivial solutions. So we only consider the equations (1.4-1.6)

Suppose that $(x,y,z,m,n)$ is a nontrivial solution of equations (1.4) or (1.5), It is clear that $m>n.$ For the equation (1.6), we have
$$x^{m}+y^{m}\geq x^{n}-y^{n}=(x-y)(x^{n-1}+x^{n-2}y+\cdots+xy^{n-2}+y^{^{n-1}}).$$ Then $x^{m}+y^{m}\geq x^{n-1}+y^{n-1}.$ So $m\geq n-1.$

If $m=n-1$ then $x-y=1, n=2, m=1.$ So the equation (1.6) become $\varphi(z)=z.$ Then the equation (1.6) has solution $(x,y,z,m,n)=(a+1,a,1,1,2).$ Correspondingly, the equation (1.2) has a solution $(a\pm1, -a, 1, 2)$ with $a\geq 1.$

We always assume that $x>y\geq1, m>n\geq 1$ in follow section.
\begin{thm}
let $\beta\geq 1, s\geq 0, a,\beta,s$ are integers. Then

(1)The only nontrivial slutions of the equation (1.4) in positive integers $x,y,z,m,n$ are $(x,y,z,m,n)=(2,1,3^{s}2^{\beta},3,1).$

(2)The only nontrivial slutions of the equation (1.5) in positive integers $x,y,z,m,n$ are $(x,y,z,m,n)=(a+1, a, 1, 2, 1), (a+2,a,2^{s},2, 1), (a+3,a,2^{\beta}3^{s},2,1).$

(3)The only nontrivial slutions of the equation (1.6) in positive integers $x,y,z,m,n$ with $\nu_{2}(x)\neq \nu_{2}(y)$ are $(x,y,z,m,n)=(2,1,p^{s}2^{\beta},q,q-1),$ where $q,p=\frac{2^{q}+1}{3}$ are both primes.
\end{thm}
By the Theorem 1.3(2), the equation (1.2) has the solution in integers $(x,y,m,n)=(a\pm i, -a, 2, 1)$ where $a$ is a integer with $i=1,2.$

The equation (1.1) has the solution in integers $(x,y,m,n)=((2^{t-1}\pm1),-(2^{t-1}\mp1),2,1), (-(2^{t-1}\pm1),(2^{t-1}\mp1),2,1)$ with $t\geq 2.$

\begin{thm}
The equations (1.6) has no nontrivial solutions in positive integers $x,y,m,n$ with $1\leq z\leq x+y, z\neq2.$
\end{thm}

Theorem 1.1 and 1.2 follow from Theorem 1.3 and Theorem 1.4 by taking $z=x+y$ and $z=1,$ respectively.

We always use the equation (1.7) to represent the equations (1.4-1.6).
$$\varphi(z\frac{x^{m}\pm y^{m}}{x+y})=z\frac{x^{n}\pm y^{n}}{x+y} \eqno{(1.7)}$$

Then we reduce Theorem 1.3 and Theorem 1.4 to the case gcd$(m,n)=1.$

Let gcd$(m,n)=d_{0}, m=d_{0}m_{0}, n=d_{0}n_{0}, x^{d_{0}}=x_{0}, y^{d_{0}}=y_{0}.$ It follow from the equation (1.4-1.6) that $2\nmid d_{0}.$
So $$z\frac{x^{d_{0}}+y^{d_{0}}}{x+y}=z_{0}.$$

Then gcd$(m_{0},n_{0})=1, m_{0}>n_{0}\geq1.$
$$\varphi(z_{0}\frac{x_{0}^{m_{0}}\pm y_{0}^{m_{0}}}{x_{0}+y_{0}})=z_{0}\frac{x_{0}^{n_{0}}\pm y_{0}^{n_{0}}}{x_{0}+y_{0}}.$$

For Theorem 1.3, noting that $\nu_{2}(x_{0})=d_{0}\nu_{2}(x)\neq d_{0}\nu_{2}(y)=\nu_{2}(y_{0})$ is equivalent to $\nu_{2}(x)\neq \nu_{2}(y).$

i)Suppose that the only nontrivial solutions of (1.4) in positive integers $x_{0},y_{0},z_{0},m_{0},n_{0}$ with gcd$(m_{0},n_{0})=1$ are $(x_{0},y_{0},z_{0},m_{0},n_{0})=(2,1,3^{s}2^{t},3,1),$
then gcd$(x_{0},y_{0})=1.$ gcd$(x_{0},y_{0})=$gcd$(x^{d_{0}},x^{d_{0}})=$gcd$(x,y)^{d_{0}}=1.$
So $d_{0}=1.$ Thus the only nontrivial solutions of the equation (1.3) in positive integers $x,y,z,m,n$ are $(x,y,z,m,n)=(2,1,3^{s}2^{t},3,1).$ ii),iii)Similar to i).

Suppose the Theorem 1.4 is true when gcd$(m, n)=1.$ and $(x,y,z,m,n)$ is a nontrivial solution of the equation(1.6) in positive integers $x,y,z,m,n$ with $1\leq z\leq x+y, z\neq2,$ gcd$(m, n)=d_{0}.$
Then $$1\leq z_{0}=z\frac{x^{d_{0}}+y^{d_{0}}}{x+y}\leq x^{d_{0}}+y^{d_{0}}=x_{0}+y_{0}, z_{0}\neq2.$$

Thus $(x_{0},y_{0},z_{0},m_{0},n_{0})$ is a nontrivial solution of the equation(1.6) in positive integers $x_{0},y_{0},z_{0},m_{0},n_{0}$ with $1\leq z_{0}\leq x_{0}+y_{0}$ and gcd$(m_{0},n_{0})=1,$ a contradiction

Let gcd$(x,y)=d_{1}, x=x_{1}d_{1}, y=y_{1}d_{1}.$ For each prime $p\geq3, p\nmid x_{1}y_{1},$ $$l_{p}=\min\{l:p|x_{1}^{l}-y_{1}^{l}\}.$$ Then $p|x_{1}^{m}-y_{1}^{m}$ is equivalent to $l_{p}\mid m.$ In addition, by Fermat theorem, $$p|x_{1}^{p-1}-y_{1}^{p-1},$$ So $l_{p}|p-1.$

Let $p, q, \gamma$ be primes. Let $p(m)$ be the least prime divisor of $m,$ $\tau(m)$ be the number of positive divisors of $m.$

Let $z=q^{\beta}z_{q}, d_{1}=q^{\alpha}d_{q}.$ then $q\nmid z_{q}d_{q}$ equation $(1.7)$ becomes
$$\varphi(q^{\beta+(m-1)\alpha}z_{q}d_{q}^{m-1}\frac{x_{1}^{m}\pm y_{1}^{m}}{x_{1}+y_{1}})=q^{\beta+(n-1)\alpha}z_{q}d_{q}^{n-1}\frac{x_{1}^{n}\pm y_{1}^{n}}{x_{1}+y_{1}} \eqno{(1.8)}$$

We prove Theorem 1.3 in section 2. Then we give some lemmas in section 3. We prove Theorem 1.4 in section 4.

\section{Proof of Theorem 1.3}
$ $

Suppose that $(x,y,z,m,n)$ is a nontrivial solution of the equation (1.7) in positive integers $x,y,z,m,n,$ then $x>y\geq 1, m>n\geq 1,$ gcd$(m,n)=1,$ So $x_{1}>y_{1}\geq 1.$ Noting that $$\frac{x_{1}^{m}+y_{1}^{m}}{x_{1}+y_{1}}>1.$$
We suppose that $$\frac{x_{1}^{m}-y_{1}^{m}}{x_{1}+y_{1}}>1.$$ ($\frac{x_{1}^{m}-y_{1}^{m}}{x_{1}+y_{1}}=1$ see Remark 2.2)

Let us discuss the situation of $\nu_{p}(x)\neq\nu_{p}(y)$ firstly. It follow that
$$2\nmid\frac{x_{1}^{m}\pm y_{1}^{m}}{x_{1}+y_{1}}, \quad 2\nmid\frac{x_{1}^{n}\pm y_{1}^{n}}{x_{1}+y_{1}}$$

By formula (1.8), let $$A=z_{2}d_{2}^{m-1}\frac{x_{1}^{m}\pm y_{1}^{m}}{x_{1}+y_{1}}, B=z_{2}d_{2}^{n-1}\frac{x_{1}^{n}\pm y_{1}^{n}}{x_{1}+y_{1}},$$

then $$\varphi(2^{\beta+\alpha(m-1)}A)=2^{\beta+\alpha(n-1)}B, 2\nmid AB, A>1, \eqno{(2.1)}.$$

If $\alpha=\beta=0,$ then (2.1) becomes $\varphi(A)=B.$ since $2\nmid AB,$ it follows that $A=B=1,$ a contradiction.

If $\alpha+\beta\geq 1.$ since $2|\varphi(A),2\nmid B$ Then (2.1) becomes $$2^{\beta+\alpha(m-1)}\frac{\varphi(A)}{2}=2^{\beta+\alpha(n-1)}B,$$
Noting that $2\nmid AB, m>n,$ we have $\alpha=0.$ Thus $\varphi(A)=2B.$ Hence there exist an prime $p\neq 2,$ a positive integer $t$ and nonnegative integers $\mu,\nu,\kappa$ such that
$$A=p^{t}=z_{2}d_{2}^{m-1}\frac{x_{1}^{m}\pm y_{1}^{m}}{x_{1}+y_{1}}, 2B=(p-1)p^{t-1}=z_{2}d_{2}^{n-1}\frac{x_{1}^{n}\pm y_{1}^{n}}{x_{1}+y_{1}}.$$
$$z_{2}=p^{\mu}, d_{2}=p^{\nu}, \frac{x_{1}^{m}\pm y_{1}^{m}}{x_{1}+y_{1}}=p^{\kappa}, t=\mu+(m-1)\nu+\kappa,$$
$$p^{t-1-\mu-(n-1)\nu}(p-1)=2\frac{x_{1}^{n}\pm y_{1}^{n}}{x_{1}+y_{1}}.$$
By $$\frac{x_{1}^{m}\pm y_{1}^{m}}{x_{1}+y_{1}}>1,$$ we have $\kappa\geq 1$ and $t-1-\mu-(n-1)\nu\geq 0.$
So $$p|\frac{x_{1}^{m}\pm y_{1}^{m}}{x_{1}+y_{1}}.$$

(1) For the equation (1.4), if $m=3,$ then $n=1,$ that is $p^{t-1-\mu-(n-1)\nu}(p-1)=2\frac{x_{1}+y_{1}}{x_{1}+y_{1}}=2.$ We have $t-1-\mu-(n-1)\nu=0.$

If $m\geq 5,$ by Carmichael primitive divisor theorem[6], we have $t-1-\mu-(n-1)\nu=0.$

In any way $t-1-\mu-(n-1)\nu=\mu+(m-1)\nu+\kappa-1-\mu-(n-1)\nu=(m-n)\nu+\kappa-1=0,$ that is $\kappa=1, \nu=0.$ Thus $d_{1}=2^{\alpha}d_{2}=1.$
$$p=\frac{x_{1}^{m}+y_{1}^{m}}{x_{1}+y_{1}}=\frac{x^{m}+y^{m}}{x+y},$$
\quad $$p-1=2\frac{x_{1}^{n}+y_{1}^{n}}{x_{1}+y_{1}}=2\frac{x^{n}+y^{n}}{x+y}$$

We have $$p-1=\frac{x^{m}+y^{m}}{x+y}-1\geq \frac{x^{n+2}+y^{n+2}}{x+y}-1\geq 2\frac{x^{n}+y^{n}}{x+y}=p-1.$$

If $x>y\geq 2,$ it follow that $$(p-1)(x+y)=x^{n+2}+y^{n+2}-x-y>4x^{n}+4y^{n}-x-y\geq 3x^{n}+3y^{n}>2x^{n}+2y^{n}=(p-1)(x+y),$$  a contradiction. So $y=1.$

If $y=1, x\geq 3,$ it follow that $$(p-1)(x+y)\geq x^{n+2}+y^{n+2}-x-y=x^{n+2}-x\geq 9x^{n}-x\geq 8x^{n}>2x^{n}+2=(p-1)(x+y),$$ a contradiction.

So, $x=2,y=1$ it follow that $2^{n+2}-2=2(2^{n}+1),$ that is $m=3, n=1,$ the equation (1.4) become $\varphi(3z)=z.$ So the equation (1.4) has the solutions $(x,y,z,m,n)=(2,1,2^{\beta}3^{s},3,1),$ where $\beta\geq 1, s\geq 0.$

(2) For the equation (1.5), if $t-1-\mu-(n-1)\nu\geq 1,$ then $p|\frac{x_{1}^{n}+y_{1}^{n}}{x_{1}+y_{1}},$
Thus $$p|x_{1}^{(m, 2n)}-y_{1}^{(m, 2n)}=(x_{1}-y_{1})(x_{1}+y_{1}), p\nmid x_{1}y_{1}.$$

If $p|x_{1}-y_{1},$ then $$p^{t-1-\mu-(n-1)\nu}(p-1)=2\frac{x_{1}^{n}+y_{1}^{n}}{x_{1}+y_{1}}\equiv 2y_{1}^{n-1}(mod p).$$  We have$p\mid y_{1},$ a contradiction.

If $p|x_{1}+y_{1}.$ By gcd$(x_{1},y_{1})=1,$ then $p\nmid x_{1}-y_{1},$
$$p^{\kappa}=\frac{x_{1}^{m}-y_{1}^{m}}{x_{1}+y_{1}}\equiv(x_{1}-y_{1})\frac{x_{1}^{m}-y_{1}^{m}}{x_{1}^{2}+y_{1}^{2}}
\equiv(x_{1}-y_{1})((x_{1}^{2})^{\frac{m}{2}-1}+(x_{1}^{2})^{\frac{m}{2}-2}y_{1}^{2}+\cdots+(y_{1}^{2})^{\frac{m}{2}-1}) \equiv(x_{1}-y_{1})\frac{m}{2}y_{1}^{m-2}(mod p), $$

$$p^{t-1-\mu-(n-1)\nu}(p-1)=2\frac{x_{1}^{n}+y_{1}^{n}}{x_{1}+y_{1}}\equiv2(x_{1}^{n-1}-x_{1}^{n-2}y_{1}+\cdots+y_{1}^{n-1})\equiv 2ny_{1}^{n-1}(mod p).$$
It follow that $p|m, p|n,$ that is $p|(m, n)=1,$ a contradiction.

So $$0=t-1-\mu-(n-1)\nu=\mu+(m-1)\nu+\kappa-1-\mu-(n-1)\nu=(m-1)\nu+\kappa-1.$$  That is $\kappa=1, \nu=0.$ Thus $d_{1}=2^{\alpha}d_{2}=1,$

$$p=\frac{x_{1}^{m}-y_{1}^{m}}{x_{1}+y_{1}}=\frac{x^{m}-y^{m}}{x+y}=(x-y)\frac{x^{m}-y^{m}}{x^{2}-y^{2}}$$
\quad $$p-1=2\frac{x_{1}^{n}+y_{1}^{n}}{x_{1}+y_{1}}=2\frac{x^{n}+y^{n}}{x+y}$$

If $x-y=1,$ we have $$p=\frac{(y+1)^{m}-y^{m}}{(y+1)^{2}-y^{2}}, p-1=2\frac{(y+1)^{n}+y^{n}}{(y+1)+y}$$ and $2\|p-1.$

i) $2|y,$ By $$p-1=\frac{(y+1)^{m}-y^{m}}{(y+1)^{2}-y^{2}}-1=\frac{(y+1)^{m}-1-y^{m}-2y}{(y+1)^{2}-y^{2}},$$ we have $4|p-1,$ a contradiction.

ii) $2\nmid y,$ By $$p-1=\frac{(y+1)^{m}-y^{m}}{(y+1)^{2}-y^{2}}-1=\frac{(y+1)^{m}-(y^{m}-1)-2(y+1)}{(y+1)^{2}-y^{2}},$$ we have $4|p-1,$ a contradiction.

So $\frac{x^{m}-y^{m}}{x^{2}-y^{2}}=1, $ then $m=2, n=1, p=x-y=3.$ The equation (1.5) becomes $\varphi(3z)=z.$  So $z=2^{\beta}3^{s}, \beta\geq 1, s\geq 0.$ The equation (1.5) has the solutions $(3+a,a,2^{\beta}3^{s},2,1)$ with $a\geq1, \beta\geq 1, s\geq 0.$

(3) For the equation (1.6), if $m=3,$ then $n=2.$
$$p^{t-1-\mu-(n-1)\nu}(p-1)=2\frac{x_{1}^{2}-y_{1}^{2}}{x_{1}+y_{1}}=2(x_{1}-y_{1}).$$ If $p\mid x_{1}-y_{1}.$ By gcd$(x_{1},y_{1})=1, p\nmid x_{1}y_{1},$ $$p^{\kappa}=\frac{x_{1}^{3}+y_{1}^{3}}{x_{1}+y_{1}}=x_{1}^{2}-x_{1}y_{1}+y_{1}^{2}\equiv y_{1}^{2},$$ a contradiction.

So $t-1-\mu-(n-1)\nu=0.$

If $m\geq 5,$ by Carmichael primitive divisor theorem, we have $t-1-\mu-(n-1)\nu=0.$

In any way $0=t-1-\mu-(n-1)\nu=\mu+(m-1)\nu+\kappa-1-\mu-(n-1)\nu=(m-n)\nu+\kappa-1,$ then $\kappa=1, \nu=0.$ Thus $d_{1}=2^{\alpha}d_{2}=1,$
$$p=\frac{x_{1}^{m}+y_{1}^{m}}{x_{1}+y_{1}}=\frac{x^{m}+y^{m}}{x+y}$$ $$p-1=2\frac{x_{1}^{n}-y_{1}^{n}}{x_{1}+y_{1}}=2\frac{x^{n}-y^{n}}{x+y},$$

if $x>y\geq 2,$ it follow that $$(x+y)(p-1)=x^{m}+y^{m}-x-y\geq 3x^{n}+2y^{n}-x-y\geq 2x^{n}+y^{n}>2x^{n}-2y^{n}=(x+y)(p-1),$$ a contradiction, So $y=1.$

if $y=1, x\geq 3,$ it follow that $$(p-1)(x+y)=x^{m}+y^{m}-x-y\geq x^{n+1}+y^{n+1}-x-y\geq3x^{n}-x>2x^{n}-2=(p-1)(x+y)$$, a contradiction.

So $x=2, y=1, n=q-1, m=q, p=\frac{2^{m}+1}{3}$ is a prime. it follow that $\varphi(pz)=(p-1)\frac{z}{2},$ we have $(x,y,z,m,n)=(2,1,2^{\beta}p^{s},q,q-1),$ where $\beta\geq 1, s\geq 0, p=\frac{2^{q}+1}{3}$ is a prime.

\begin{rem}
For the equation (1.4), we have $2\nmid mn,$  if $\nu_{2}(x)=\nu_{2}(y),$ then $2\nmid x_{1}y_{1},$ we have $$2\nmid\frac{x_{1}^{m}+y_{1}^{m}}{x_{1}+y_{1}}, 2\nmid\frac{x_{1}^{n}+y_{1}^{n}}{x_{1}+y_{1}}.$$
So the above discussion (1) of the equation (1.4) can remove the condition $\nu_{2}(x)\neq\nu_{2}(y).$
\end{rem}

\begin{rem}
If$$\frac{x_{1}^{m}-y_{1}^{m}}{x_{1}+y_{1}}=1,$$ that is the equation (1.5).
then $m=2, n=1, x_{1}-y_{1}=1.$
By fomula (1.8), the equation (1.5) becomes $$ \varphi(2^{\beta+\alpha}z_{2}d_{2})=2^{\beta}z_{2}.$$

i)$\alpha=\beta=0,$ then $\varphi(z_{2}d_{2})=z_{2}.$ So $z=z_{2}=1, d_{1}=d_{2}=1,$ the equation (1.5) becomes $\varphi(z)=z.$ So the equation (1.5) has the integer solution $(a+1, a, 1, 2, 1)$ with $a\geq 1.$

ii) $\alpha+\beta>0,$ then $2^{\alpha-1}\varphi(z_{2}d_{2})=z_{2},$ So $\alpha\leq 1.$

If $\alpha=1,$ then $\varphi(z_{2}d_{2})=z_{2},$ $z=2^{\beta}z_{2}=2^{\beta}, d_{1}=2d_{2}=2, \varphi(2z)=z.$
The equation (1.5) has the integer solutions $(2a+2, 2a, 2^{s}, 2, 1),$ where $a\geq 1, s\geq 0.$

If $\alpha=0, \beta\geq 1$ then $\varphi(z_{2}d_{2})=2z_{2},$ $z=2^{\beta}z_{2}=2^{\beta}3^{s}, d_{1}=d_{2}=3.$
The equation (1.5) has the integer solutions $(3a+3, 3a, 2^{\beta}3^{s}, 2, 1),$ where $a\geq 1, \beta\geq 1, s\geq 0.$

For the equation (1.5), we have $2|m, 2\nmid n,$ if $\nu_{2}(x)=\nu_{2}(y),$ then $2\nmid x_{1}y_{1},$
$$2|\frac{x_{1}^{m}-y_{1}^{m}}{x_{1}+y_{1}}=(x_{1}-y_{1})\frac{x_{1}^{m}-y_{1}^{m}}{x_{1}^{2}-y_{1}^{2}}, 2\nmid\frac{x_{1}^{n}+y_{1}^{n}}{x_{1}+y_{1}},$$

It follow that $2|A, 2\nmid B,$ since $\varphi(2^{\beta+\alpha(m-1)+1}\frac{A}{2})=2^{\beta+\alpha(n-1)}B,$ then $2\parallel A,$ $2^{\alpha(m-n)}\varphi(\frac{A}{2})=B.$ So $\alpha=0, A=2, B=1, \varphi(2z)=z$ We have $(x,y,z,m,n)=(2a+1,2a-1,2^{s},2,1), a\geq 1, s\geq 0.$
\end{rem}

This completes the proof of Theorem 1.3. \hfill$\Box$\\

\section{some Lemmas}
$ $

In order to prove Theorem 1.4, we give some lemmas in this section. We always assume that $(x,y,z,m,n)$ is a non trivial solution of the equation (1.6)$$\varphi(z\frac{x^{m}+y^{m}}{x+y})=z\frac{x^{n}-y^{n}}{x+y}, 2\nmid m, 2|n,$$ with $1\leq z\leq x+y, z\neq2.$ It follow from Theorem 1.3(3) that Theorem 1.4 is true when $\nu_{p}(x)\neq\nu_{p}(y).$ So we always assume that $\nu_{2}(x)=\nu_{2}(y),$ then $x_{1}, y_{1}$ are both odd.

Corresponding to the lemmas of Equation $$\varphi(z\frac{x^{m}-y^{m}}{x-y})=z\frac{x^{n}-y^{n}}{x-y},$$ with $1\leq z\leq x-y,$
in [8], we get some similar lemmas about the equation (1.6).

By $m>n,$ we have
$$x<\frac{z(x^{m}+y^{m})/(x+y)}{z(x^{n}-y^{n})/(x+y)}
=\frac{z(x^{m}+y^{m})/(x+y)}{\varphi(z(x^{m}+y^{m})/(x+y))}=\prod_{p|z\frac{x^{m}+y^{m}}{x+y}}(1+\frac{1}{p-1}) \eqno{(3.1)}$$

\begin{lemma}\label{lemma1}
Let q be a prime and $q|m,$ then

(1)$$q^{\frac{1}{2}\nu_{q}(m)\tau(m)-1}|x_{1}^{q-1}-y_{1}^{q-1}.$$
If $q\nmid z,$ then $$q^{\frac{1}{2}\nu_{q}(m)\tau(m)}|x_{1}^{q-1}-y_{1}^{q-1}.$$

(2)If $q|x_{1}-y_{1},$ then $\frac{1}{2}\nu_{q}(m)\tau(m)-1\leq\nu_{q}(x_{1}-y_{1}).$

(3)If $k$ is a positive integer such that $$q^{k}\nmid x_{1}^{q-1}-y_{1}^{q-1},$$ then there are at most $k$ distinct primes $p$ with $l_{p}|2m,l_{p}\nmid m, q|l_{p}.$
\end{lemma}

Proof.
(1) Let $m=q^{\nu_{q}(m)}m_{q}$ and $l_{1}, l_{2},\cdots, l_{t}$ be all positive divisors of $m_{q}.$ Then $q^{i}l_{j}(1\leq i\leq \nu_{q}(m),1\leq j\leq t)$ are all distinct positive divisors of $m.$ By Carmichael primitive divisor theorem(see[6]), each of $x_{1}^{q^{i}l_{j}}+y_{1}^{q^{i}l_{j}}$ has a primitive prime divisor $p_{i,j}\equiv 1(mod q^{i}l_{j}).$ It is clear that
$$p_{i,j}|\frac{x_{1}^{q^{i}l_{j}}+y_{1}^{q^{i}l_{j}}}{x_{1}+y_{1}}, \qquad \frac{x_{1}^{q^{i}l_{j}}+y_{1}^{q^{i}l_{j}}}{x_{1}+y_{1}}| \frac{x_{1}^{m}+y_{1}^{m}}{x_{1}+y_{1}}.$$

By formla (1.8), we have $$q^{\beta+(m-1)\alpha}\prod_{1\leq i\leq \nu_{q}(m), 1\leq j\leq t}p_{i,j}|q^{\beta+(m-1)\alpha}z_{q}d_{q}^{m-1}\frac{x_{1}^{m}+y_{1}^{m}}{x_{1}+y_{1}}.$$

So $$\varphi(q^{\beta+(m-1)\alpha})\prod_{1\leq i\leq \nu_{q}(m), 1\leq j\leq t}(p_{i,j}-1)|\varphi(q^{\beta+(n-1)\alpha}z_{q}d_{q}^{n-1}\frac{x^{m}+y^{m}}{x+y}).$$

By formla (1.8), $$\varphi(q^{\nu_{q}(z)+(m-1)\alpha})\prod_{1\leq i\leq \nu_{q}(m), 1\leq j\leq t}q^{i}|q^{\beta+(m-1)\alpha}\frac{x_{1}^{n}-y_{1}^{n}}{x_{1}-y_{1}}.$$

We have $$\varphi(q^{\beta+(m-1)\alpha})q^{\frac{1}{2}\nu_{q}(m)\tau(m)}|q^{\beta+(n-1)\alpha}\frac{x_{1}^{n}-y_{1}^{n}}{x_{1}-y_{1}}\eqno{(3.2)}.$$

We divide into two cases.

Case 1: $\alpha\geq 1$ or $\alpha=\beta=0,$ we have $$q^{\frac{1}{2}\nu_{q}(m)\tau(m)}|(x_{1}^{n}-y_{1}^{n}).$$

Since gcd$(x_{1}, y_{1})=1,$ it follow that $q\nmid x_{1}y_{1}.$ By Euler theorem,
$$q^{\frac{1}{2}\nu_{q}(m)\tau(m)}|x_{1}^{(q-1)\frac{1}{2}\nu_{q}(m)\tau(m)-1}-y_{1}^{(q-1)\frac{1}{2}\nu_{q}(m)\tau(m)-1}.$$

So
$$q^{\frac{1}{2}\nu_{q}(m)\tau(m)}|x_{1}^{((q-1)\frac{1}{2}\nu_{q}(m)\tau(m)-1,n)}-y_{1}^{((q-1)\frac{1}{2}\nu_{q}(m)\tau(m)-1,n)}.$$
By gcd$(m,n)=1,$ we have gcd$(q,n)=1.$ It follow that
$$q^{\frac{1}{2}\nu_{q}(m)\tau(m)}|x_{1}^{(q-1,n)}-y_{1}^{(q-1,n)},$$ That is $$q^{\frac{1}{2}\nu_{q}(m)\tau(m)}|x_{1}^{q-1}-y_{1}^{q-1}.$$

Case 2: $\alpha=0,\beta\geq 1.$
So we have $$q^{\beta-1}q^{\frac{1}{2}\nu_{q}(m)\tau(m)}|q^{\beta}(x_{1}^{n}-y_{1}^{n}).$$
It follow that $$q^{\frac{1}{2}\nu_{q}(m)\tau(m)-1}|x_{1}^{n}-y_{1}^{n}.$$
If $\frac{1}{2}\nu_{q}(m)\tau(m)-1 \geq 1,$ then, similar to Case 1, we have $$q^{\frac{1}{2}\nu_{q}(m)\tau(m)-1}|x_{1}^{q-1}-y_{1}^{q-1}.$$
It is clear that it also holds if $\frac{1}{2}\nu_{q}(m)\tau(m)-1=0.$

(2)If $q|x_{1}-y_{1},$ noting that $2|n,$ then
\begin{eqnarray*}
\frac{x_{1}^{n}-y_{1}^{n}}{x_{1}+y_{1}} &=& (x_{1}-y_{1})\frac{x_{1}^{n}-y_{1}^{n}}{x_{1}^{2}-y_{1}}^{2}\\
&=&(x_{1}-y_{1})((x_{1}^{2})^{\frac{n}{2}-1}+(x_{1}^{2})^{\frac{n}{2}-2}y_{1}^{2}+\cdots+x_{1}^{2}(y_{1}^{2})^{\frac{n}{2}-2}y_{1}^{2}+(x_{1}^{2})^{\frac{n}{2}-1})\\
&\equiv &(x_{1}-y_{1})\frac{n}{2}(x_{1}^{2})^{\frac{n}{2}-1}(mod q)\\
\end{eqnarray*}

Since $q\nmid x_{1}y_{1}, q\nmid n,$ it follow that $$\nu_{q}(\frac{x_{1}^{n}-y_{1}^{n}}{x_{1}+y_{1}} )=\nu_{q}(x_{1}-y_{1}).$$

By(3.2), we have
$$\varphi(q^{\beta+(m-1)\alpha})q^{\frac{1}{2}\nu_{q}(m)\tau(m)}|q^{\beta+(n-1)\alpha}\nu_{q}(x_{1}-y_{1}).$$

That is $$\frac{1}{2}\nu_{q}(m)\tau(m)-1\leq \nu_{q}(x_{1}-y_{1}).$$

(3)Suppose that there are at least $k+1$ primes $p$ with $l_{p}|2m,l_{p}\nmid m, q|l_{p}.$ Let $p_{1}, p_{2}, \cdots, p_{k+1}$ be $k+1$ distibct primes with $l_{p_{i}}|2m,l_{p_{i}}\nmid m, q|l_{p_{i}}.$ Then, for $1\leq i\leq k+1,$

$$p_{i}|\frac{x_{1}^{l_{p_{i}}}+y_{1}^{l_{p_{i}}}}{x_{1}+y_{1}},\qquad \frac{x_{1}^{l_{p_{i}}}+x_{1}^{l_{p_{i}}}}{x_{1}+y_{1}}|\frac{x_{1}^{m}+y_{1}^{m}}{x_{1}+y_{1}}.$$

By (1.8), $$q^{\beta+(m-1)\alpha}p_{1}p_{2}\cdots p_{k+1}|zd_{1}^{m-1}\frac{x_{1}^{m}+y_{1}^{m}}{x_{1}+y_{1}}$$
Then
$$\varphi(q^{\beta+(m-1)\alpha})(p_{1}-1)(p_{2}-1)\cdots(p_{k+1}-1)|\varphi(zd_{1}^{m-1}\frac{x_{1}^{m}+y_{1}^{m}}{x_{1}+y_{1}}).$$

Noting $p_{i}\nmid x_{1}y_{1}(1\leq i\leq k+1), q|l_{p_{i}},$ it follows that $q|p_{i}-1(1\leq i\leq k+1).$
$$q^{\beta+(m-1)\alpha+k}|\varphi(zd^{m-1}\frac{x_{1}^{m}+y_{1}^{m}}{x_{1}+y_{1}})$$

It follow from $(1.8),$ $$q^{\beta+(m-1)\alpha+k}|q^{\beta+(m-1)\alpha}\frac{x_{1}^{n}-y_{1}^{n}}{x_{1}+y_{1}}.$$

Noting that $m>n,$ It follow that $$q^{k}|x_{1}^{n}-y_{1}^{n}.$$

It follow from $q\nmid x_{1}y_{1}.$ Similar to (1)case 1, $$q^{k}|x_{1}^{q-1}-y_{1}^{q-1},$$ a contradiction.

This completes the proof of the Lemma 3.1.

\begin{lemma}\label{lemma1}
Let $\omega(m)$ be the number of prime divisors of $m,$ and let $\tau(m)$ is the number of positive divisors of $m.$
Then $$\tau(m)<2\max\{p(m), x\}$$ and $$\omega(m)<\frac{\log(2\max\{p(m), x\})}{\log 2}$$
\end{lemma}

Proof. In this proof, $p(m)$ is abbreviated as $p.$ If $p|z,$ then by $p|m$ and lemma 4.1, noting that $1\leq z\leq x+y\leq 2x,$ we have
$$p^{\frac{1}{2}\tau(m)-1}\leq p^{\frac{1}{2}\nu_{p}(m)\tau(m)-1}\leq x_{1}^{p-1}-y_{1}^{p-1}\leq x^{p-1}$$

So $$\tau(m)<\frac{2p\log x-2\log x+2\log p}{\log p} \eqno{(3.3)}.$$

Use another estimate, we have $$p^{\frac{1}{2}\tau(m)-1}\leq x^{p-1}<\frac{2}{z}x^{p}\leq2p^{-1}x^{p}.$$ So
$$\tau(m)<\frac{2p\log x+2\log 2}{\log p}  \eqno{(3.4)}.$$

If $p\nmid z,$ then by $p|m$ and Lemma 4.1(1), we have $$p^{\frac{1}{2}\nu_{p}(m)\tau(m)}|x_{1}^{p-1}-y_{1}^{p-1}.$$

Hence $$p^{\frac{1}{2}\tau(m)}\leq p^{\frac{1}{2}\nu_{p}(m)\tau(m)}\leq x_{1}^{p-1}-y_{1}^{p-1}\leq x^{p-1}.$$
So $$\tau(m)<\frac{2p\log x-2\log x}{\log p}\eqno{(3.5)}.$$

If $p\leq x,$ then by $p\geq 3,$ $$\frac{2p}{\log p}<\frac{2x}{\log x}.$$ It follow from (3.3) and (3.5)that $\tau(m)<2x.$

If $p> x,$ then $$(\frac{p}{x})^{p}\geq(1+\frac{1}{x})^{p}>1+p\frac{1}{x}\geq 2.$$ It follow from (3.4) and (3.5)that $\tau(m)<2p.$

This completes the proof of Lemma 3.2. \hfill$\Box$\\

\begin{lemma}\label{lemma1}
Let $d$ be a divisor of $2m$ with $d>30$ and let $$S_{d}=\sum_{l_{p}=d}\frac{1}{p}.$$ Then $$S_{d}<\frac{1.084}{d}+\frac{1}{d\log(d+1)}+\frac{2\log\log d}{\varphi(d)}+\frac{2\log\log x}{\varphi(d)\log d}.$$
\end{lemma}

Proof. We follow the proof of [8, lemma 4.7]. Let $P_{d}=\{p: l_{p}=d\},$ then $d=l_{p}|p-1$ for all $p\in P_{d}.$ Hence
$$(d+1)^{\sharp P_{d}}\leq\prod_{P_{d}}p\leq x_{1}^{d}-y_{1}^{d}\leq x^{d}.$$
It follow that $$\sharp P_{d}\leq\frac{d\log x}{\log(d+1)}.$$

Let $$\pi(X;d,1)=\sharp\{p:p\leq X, d\mid p-1\},$$ by the Brun-Titchmarsh theorem due to Montgomery and Vaughan[23],
$$\pi(X;d,1)<\frac{2X}{\varphi(d)\log(\frac{X}{d})},$$
for all $X>d\geq2.$

Let $A_{d}=\{p:p\leq 4d, d\mid p-1\}.$ We split $S_{d}$ as follows:
\begin{eqnarray*}
S_{d}&=& \sum_{p\leq4d,l_{p}=d}\frac{1}{p}+\sum_{4d<p\leq d^{2}\log x,l_{p}=d}\frac{1}{p}+\sum_{p>d^{2}\log x,l_{p}=d}\frac{1}{p}\\
&\leq & \sum_{p\in A_{d}}\frac{1}{p}+\sum_{4d<p\leq d^{2}\log x,d\mid p-1}\frac{1}{p}+\sum_{p>d^{2}\log x,p\in P_{d}}\frac{1}{p}\\
&\leq & T_{1}+T_{2}+T_{3}.\\
\end{eqnarray*}

For $T_{2},$ we have

\begin{eqnarray*}
T_{2}&=& \int^{d^{2}\log x}_{4d}\frac{1}{t}d\pi(X;d,1)\\
&= & \frac{\pi(X;d,1)}{t}|^{d^{2}\log x}_{4d}+\int^{d^{2}\log x}_{4d}\frac{\pi(X;d,1)}{t^{2}}dt\\
&\leq & \frac{2}{\varphi(d)\log(d\log x)}-\frac{\pi(X;d,1)}{4d}+\frac{2}{\varphi(d)}\int^{d^{2}\log x}_{4d}\frac{1}{t\log(\frac{t}{d})}dt.\\
&\leq & \frac{2\log\log(d\log x)}{\varphi(d)}-\frac{\pi(X;d,1)}{4d}+\frac{2}{\varphi(d)}(\frac{1}{\log(d\log x)}-\log\log4).\\
\end{eqnarray*}

Since $d\geq30$ and $x\geq x_{1}\geq3,$ it follow that
$$\frac{1}{\log(d\log x)}-\log\log4<\frac{1}{\log(30)}-\log\log4<0.$$
Hence
$$T_{1}+T_{2}\leq\frac{2\log\log(d\log x)}{\varphi(d)}-\frac{\pi(X;d,1)}{4d}+\sum_{p\in A_{d}}\frac{1}{p}.$$
By$A_{d}\subseteq\{d+1,2d+1,3d+1\},\pi(X;d,1)=\sharp A_{d}\leq3.$
It follow that
$$-\frac{\pi(X;d,1)}{4d}+\sum_{p\in A_{d}}\frac{1}{p}\leq -\frac{3}{4d}+\frac{1}{d+1}+\frac{1}{2d+1}+\frac{1}{3d+1}<\frac{1.084}{d}.$$
So $$T_{1}+T_{2}\leq\frac{1.084}{d}+\frac{2\log\log(d\log x)}{\varphi(d)}.$$
For $T_{3},$ by(4.6), $$T_{3}<\frac{\sharp P_{d}}{d^{2}\log x}<\frac{1}{d\log(d+1)}.$$

Therefore,$$T_{1}+T_{2}+T_{3}<\frac{1.084}{d}+\frac{1}{d\log(d+1)}+\frac{2\log\log(d\log x)}{\varphi(d)}.$$

Noting that
\begin{eqnarray*}
\log\log(d\log x)&=& \log(\log d+\log\log x)\\
&=& \log\log d+\log(1+\frac{\log\log x}{\log d})\\
&<& \log\log d+\frac{\log\log x}{\log d}\\
\end{eqnarray*}
we have $$S_{d}<\frac{1.084}{d}+\frac{1}{d\log(d+1)}+\frac{2\log\log d}{\varphi(d)}+\frac{2\log\log x}{\varphi(d)\log d}.$$

This completes the proof of Lemma 3.3. \hfill$\Box$\\

\begin{lemma}\label{lemma1}
We have $$x_{1}\frac{\varphi(zd_{1})}{z}<\prod_{l_{p}|2m,l_{p}\nmid m,l_{p}>2}(1+\frac{1}{p-1}).$$
\end{lemma}

Proof. By(3.1), we have
\begin{eqnarray*}
x &<& \prod_{p|z\frac{x^{m}+y^{m}}{x+y}}(1+\frac{1}{p-1})\\
&\leq &\prod_{p|zd_{1}}(1+\frac{1}{p-1})\prod_{p|\frac{x_{1}^{m}+y_{1}^{m}}{x_{1}+y_{1}}}(1+\frac{1}{p-1})\\
&= &\frac{zd_{1}}{\varphi(zd_{1})}\prod_{p|\frac{x_{1}^{m}+y_{1}^{m}}{x_{1}+y_{1}}}(1+\frac{1}{p-1})\\
\end{eqnarray*}

By[4] gcd$(\frac{x_{1}^{m}+y_{1}^{m}}{x_{1}+y_{1}},\frac{x_{1}^{m}-y_{1}^{m}}{x_{1}-y_{1}})=1,$
We have
$$x_{1}\frac{\varphi(zd_{1})}{z}<\prod_{l_{p}|2m,l_{p}\nmid m}(1+\frac{1}{p-1})<\prod_{p|\frac{x_{1}^{m}+y_{1}^{m}}{x_{1}+y_{1}}}(1+\frac{1}{p-1}).$$

If $l_{p}=1,$ then $p|x_{1}-y_{1}.$ By gcd$(x_{1},y_{1})=1,$  $p\nmid x_{1}y_{1}.$ We have $$\frac{x_{1}^{m}+y_{1}^{m}}{x_{1}+y_{1}}=x^{m-1}-x^{m-2}y_{1}+\cdots+x_{1}^{2}y_{1}^{m-3}-x_{1}y_{1}^{m-2}+y_{1}^{m-1}\equiv y_{1}^{m-1}(mod p).$$
Since $p|\frac{x_{1}^{m}+y_{1}^{m}}{x_{1}+y_{1}},$ it follow that $p|y_{1}^{m-1},$ a contradiction.

If $l_{p}=2,$ then $p\nmid x_{1}-y_{1}, p|x_{1}^{2}-y_{1}^{2}.$ By gcd$(x_{1},y_{1})=1$ We have $p|x_{1}+y_{1}, p\nmid x_{1}y_{1}.$

Then $$\frac{x_{1}^{m}+y_{1}^{m}}{x_{1}+y_{1}}=x^{m-1}-x^{m-2}y_{1}+\cdots+x_{1}^{2}y_{1}^{m-3}-x_{1}y_{1}^{m-2}+y_{1}^{m-1}\equiv my_{1}^{m-1}(mod p).$$

Since $p\nmid y_{1},$ it follow that $p|m.$ Let $m=p^{s}m_{1},$ By Carmichael primitive divisor theorem, $x_{1}^{pm_{1}}+y_{1}^{pm_{1}}$ has a primitive prime divisor $\gamma\equiv1(mod pm_{1}).$

By formula (1.8) and $p|\frac{x_{1}^{m}+y_{1}^{m}}{x_{1}+y_{1}},$ we have $$p^{\beta+(m-1)\alpha+1}\gamma|z\frac{x^{m}+y^{m}}{x+y},$$ $$\varphi(p^{\beta+(m-1)\alpha+1})(\gamma-1)|p^{\beta+(n-1)\alpha}z_{p}d_{p}\frac{x_{1}^{n}-y_{1}^{m}}{x_{1}+y_{1}}.$$

It follow that $$p^{\beta+(m-1)\alpha+1}|p^{\beta+(n-1)\alpha}\frac{x_{1}^{n}-y_{1}^{n}}{x_{1}+y_{1}}.$$

By $p|m, (m,n)=1,$ we have $p\nmid n.$ So
$$\frac{x_{1}^{n}-y_{1}^{n}}{x_{1}+y_{1}}=(x_{1}-y_{1})((x^{2})^{\frac{n}{2}-1}+(x^{2})^{\frac{n}{2}-2}y_{1}^{2}+\cdots+x_{1}^{2}(y_{1}^{2})^{\frac{n}{2}-2}+(y_{1}^{2})^{\frac{n}{2}-1})\equiv (x_{1}-y_{1})\frac{n}{2}(y_{1}^{2})^{\frac{n}{2}-1}(mod p)$$
It follow that $p\nmid \frac{x_{1}^{n}-y_{1}^{n}}{x_{1}+y_{1}}.$ So, $$p^{\beta+(m-1)\alpha+1}|p^{\beta+(n-1)\alpha}.$$ a contradiction.

This completes the proof of Lemma 3.4. \hfill$\Box$\\

\begin{lemma}\label{lemma1}
If$x\leq73$ and $d\mid 2m,$ then

(1)If $d\geq173,$  then$$\log(\prod_{l_{p}=d,d>2}(1+\frac{1}{p-1}))<3.7341\sum_{l_{p}=d,d>2}\frac{\log\log d}{\varphi(d)}.$$

(2)If $p^{\prime}(d)=\min\{p:p\neq2,p\mid d\}\geq 173,$ then $\log(\prod_{l_{p}=d, d>2}(1+\frac{1}{p-1}))<0.03834$
\end{lemma}

Proof. (1) Since $d>173$ and $x\geq x_{1}\geq3,$ it follow from lemma 3.4 that
\begin{eqnarray*}
\log(\prod_{l_{p}=d}(1+\frac{1}{p-1}))&=& \sum_{l_{p}=d}\log(1+\frac{1}{p-1})<\sum_{l_{p}=d}\frac{1}{p-1}\\
&=& \sum_{l_{p}=d}\frac{1}{p(p-1)}+\sum_{l_{p}=d}\frac{1}{p}<\sum_{n=d+1}^{\infty}\frac{1}{n(n-1)}+\sum_{l_{p}=d}\frac{1}{p}\\
&=& \frac{1}{d}+S_{d}<\frac{2.084}{d}+\frac{1}{d\log(d+1)}+\frac{2\log\log d}{\varphi(d)}+\frac{2\log\log x}{\varphi(d)\log d}.\\
&<& \frac{2.084}{\varphi(d)}+\frac{1}{\varphi(d)\log174}+\frac{2\log\log d}{\varphi(d)}+\frac{2\log\log73}{\varphi(d)\log173}.\\
&<& \frac{2.8431}{\varphi(d)}+\frac{2\log\log d}{\varphi(d)}\\
&=& \frac{3.7341\log\log d}{\varphi(d)}\\
\end{eqnarray*}

(2)let $d=2p_{1}p_{2}\cdots p_{t},$ where $2<p_{1}\leq p_{2}\leq\cdots \leq p_{t}$ are primes.
Since $p_{i}\geq p^{\prime}(d)\geq 173(1\leq i\leq t),$
it follow that $$\frac{\log\log (p_{i})}{\varphi(p_{i})}\leq\frac{\log\log173}{172}<1,(2\leq i\leq t).$$

By the lemma 4.8 of [8], If $a, b\geq 78,$ then $\log\log(ab)\leq(\log\log(a))(\log\log(b)).$
In view of (1) and lemma3.4,
\begin{eqnarray*}
\log(\prod_{l_{p}=d, d>2}(1+\frac{1}{p-1}))&\leq& \frac{3.3741\log\log d}{\varphi(d)}\leq\frac{3.3741\log\log(2p_{1})}{\varphi(2p_{1})}\cdots\frac{\log\log(p_{t})}{\varphi(p_{t})}\\
&<& \frac{3.7341\log\log (2p_{1})}{\varphi(2p_{1})}\leq\frac{3.7341\log\log (346)}{172}<0.03834.\\
\end{eqnarray*}

This completes the proof of Lemma 3.5. \hfill$\Box$\\

\begin{lemma}\label{lemma1}
If $x\leq 73, q\nmid x_{1}-y_{1},$ where $q$ is a prime with $q\nmid m, q<173,$ then $q^{3}\nmid x_{1}^{q-1}-y_{1}^{q-1}.$
\end{lemma}

Proof. We follow the proof of [8, lemma 6.9]. Since $x\leq73, 2\nmid m$ it follow that $x_{1} \leq 73, q \geq 3.$ A simple calculation by a computer shows that, for any integers $1\leq y_{1}<x_{1}\leq 73,$ there are no odd primes $3\leq p< 173$ such that $$p\nmid x_{1}-y_{1}, p^{6}| x_{1}^{q-1}-y_{1}^{q-1}.$$
So $$q^{6}\nmid x_{1}^{q-1}-y_{1}^{q-1} \eqno{(3.6)}$$

If $m$ has a prime divisor $q<173,$ by lemma 3.1(1), $$q^{\frac{1}{2}\nu_{q}(m)\tau(m)-1}| x_{1}^{q-1}-y_{1}^{q-1}.$$ Hence $\frac{1}{2}\nu_{q}(m)\tau(m)-1\leq 5.$ That is $\nu_{q}(m)\tau(m)\leq 12.$ Therefore, $\tau(m)\leq 12.$

Suppose that $ p^{3}| x_{1}^{q-1}-y_{1}^{q-1}.$ Hence $$p\nmid x_{1}-y_{1}, p^{3}| x_{1}^{q-1}-y_{1}^{q-1} \eqno{(3.7)}$$
A simple calculation by a computer shows that, for $1\leq y_{1}<x_{1}\leq 9,$ there is no prime $\gamma<173$ satisfying (3.7); for $10\leq x_{1}\leq73,$ there are at most two primes $\gamma<173$ satisfying (3.7).

It follow that $10\leq x_{1}<73.$ By lemma 3.5, we have
$$x_{1}\frac{\varphi(zd_{1})}{z}<\prod_{l_{p}|2m,l_{p}\nmid m,l_{p}>2}(1+\frac{1}{p-1})=\prod_{l_{p}\in D_{1}}(1+\frac{1}{p-1})\prod_{l_{p}\in D_{2}}(1+\frac{1}{p-1}),$$
where$$D_{1}=\{d:d|2m,d\nmid m,d>2,\exists p<173, p|d\},D_{2}=\{d:d|2m,d\nmid m,d>2,d\not\in D_{1}\}$$

It follow that, $l_{p}>2, l_{p}| p-1,$ we have $p>7, p-1\neq 2^{t}.$ By (3.7) and lemma 4.1(3), there are at most 6 primes $p$ with $l_{p}|2m,l_{p}\nmid m,l_{p}>2, q|l_{p}.$
So, for any given prime $q<173, \sharp\{p:d\in D_{1}, q|l_{p}\}\leq 6.$
Then, $ \sharp\{p:l_{p}\in D_{1}\}\leq 12.$ It follow that
$$\prod_{l_{p}\in D_{1}}(1+\frac{1}{p-1})\leq \prod_{4\leq i\leq16}^{p\neq 17}(1+\frac{1}{p-1})<1.8443,$$ where $p_{i}$ is the $i$-th prime.

By $\tau(m)\leq 12, q<173,$ we have $q\not\in D_{2}.$ Noting $1\not\in D_{2},$ we have $\sharp D_{2}\leq \tau(m)-2=10.$ So
$$\prod_{l_{p}\in D_{2}}(1+\frac{1}{p-1})=\prod_{d\in D_{2}}\exp(\log\prod_{l_{p}=d}(1+\frac{1}{p-1})\leq \prod_{d\in D_{2}}\exp0.03834\leq\exp10\times0.03834<1.4673.$$

Hence $$x_{1}\frac{\varphi(zd_{1})}{z}<\prod_{l_{p}\in D_{1}}(1+\frac{1}{p-1})\prod_{l_{p}\in D_{2}}(1+\frac{1}{p-1})<1.8443\times1.4673<2.7062.$$

It is clear that $1\leq z\leq x+y<146<2\times3\times5\times7.$ Since $10\leq x_{1}\leq 73$ and $x_{1}$ is odd, it follows that $x_{1}\geq 11.$
we have $$x_{1}\frac{\varphi(zd_{1})}{z}\geq 11\frac{\varphi(z)}{z}\geq11\frac{1}{2}\frac{2}{3}\frac{4}{5}>2.93$$ a contradiction.

This completes the proof of Lemma 3.6. \hfill$\Box$\\

\begin{lemma}\label{lemma1}
If $x\leq73,$ then $p(m)\geq173.$
\end{lemma}

Proof. We follow the proof of [8, lemma 6.10]. Suppose that $p(m)\leq173.$ Let $q$ be a prime with $q\mid m, q<173.$

If $q\nmid x_{1}-y_{1}.$ By lemma 3.6, $$q^{3}\nmid x_{1}^{q-1}-y_{1}^{q-1}.$$
By lemma 3.1(1), $$q^{\frac{1}{2}\nu_{q}(m)\tau(m)-1}|x_{1}^{q-1}-y_{1}^{q-1}.$$
Hence $\frac{1}{2}\nu_{q}(m)\tau(m)-1\leq2.$ That is, $\nu_{q}(m)\tau(m)\leq6.$ So $m\in \{q,q^{2},q\gamma,q\gamma^{2}\}.$

If $q\mid x_{1}-y_{1}.$ By lemma 4.1(2), $$\nu_{q}(x_{1}-y_{1})\geq \frac{1}{2}\nu_{q}(m)\tau(m)-1.$$
Noting $x_{1}-y_{1}\leq x-y<73,$ we have $\nu_{q}(x_{1}-y_{1})\leq3.$ So $\nu_{q}(m)\tau(m)\leq8.$
It follow that $m\in \{q,q^{2},q\gamma,q\gamma^{2},qp\gamma,q\gamma^{3}\}.$ But for $q\gamma^{3},$ if $\gamma\nmid x_{1}-y_{1},$ we have $m\neq q\gamma^{3},$ if $\gamma\mid x_{1}-y_{1},$ we have $$\frac{1}{2}\nu_{\gamma}(m)\tau(m)-1=11,$$ a contradiction.

In any way $$m\in \{q,q^{2},q\gamma,q\gamma^{2},qp\gamma\}.$$

Noting that: if $q>3,$ $\nu_{q}(x_{1}-y_{1})\leq2,$  we have  $\frac{1}{2}\nu_{q}(m)\tau(m)-1\leq2.$
It follow that $m\in \{q,q^{2},q\gamma,q\gamma^{2}\}.$

By lemma 3.5
$$x_{1}\frac{\varphi(zd_{1})}{z}<\prod_{l_{p}|2m,l_{p}\nmid m,l_{p}>2}(1+\frac{1}{p-1})=\prod_{l_{p}\in D_{1}}(1+\frac{1}{p-1})\prod_{l_{p}\in D_{2}}(1+\frac{1}{p-1}),$$
where $$D_{1}=\{d:d|2m,d\nmid m,d>2,\exists p<173, p|d\},D_{2}=\{d:d|2m,d\nmid m,d>2,d\not\in D_{1}\}.$$

Noting $2\nmid m, l_{p}|p-1$ we have $p\geq7.$ Let $p_{i}$ be the $i-$th prime, then $i\geq 4.$

Since $x_{1}$ and $y_{1}$ are odd and $x_{1}>y_{1}\geq1,$ it follow that $x_{1}\geq 3.$ So $x\geq 3.$

It is clear that $z\leq x+y\leq 146.$ We divide into three cases:

Case 1: $x>3,$ then either $x_{1}\geq5$ or $d_{1}\geq2.$ remove the case where $d_{1}=1,x_{1}=5,z=6.$
(1)If $d_{1}\geq 7$ or $d_{1}=5,$ then
$$x_{1}\frac{\varphi(zd_{1})}{z}\geq 3\varphi(d_{1})\frac{\varphi(z)}{z}\geq 3\times4\frac{1}{2}\frac{2}{3}\frac{4}{5}\geq2.8.$$
If $d=6,4,3, x_{1}\geq5$ then
$$x_{1}\frac{\varphi(zd_{1})}{z}\geq 5\varphi(d_{1})\frac{\varphi(z)}{z}\geq 5\times2\frac{1}{2}\frac{2}{3}\frac{4}{5}\geq2.$$
If $d=6,4,3, x_{1}=3$ then
$$x_{1}\frac{\varphi(zd_{1})}{z}\geq 3\varphi(d_{1})\frac{\varphi(z)}{z}\geq 3\times2\frac{1}{2}\frac{2}{3}\geq2.$$
If $d=2,1,x_{1}\geq 9,$ then $$x_{1}\frac{\varphi(zd_{1})}{z}\geq 9\frac{\varphi(z)}{z}\geq 9\frac{1}{2}\frac{2}{3}\frac{4}{5}\geq2.4.$$
If $d=2,1,x_{1}=7,$ then $$x_{1}\frac{\varphi(zd_{1})}{z}\geq 7\frac{\varphi(z)}{z}\geq 7\frac{1}{2}\frac{2}{3}\geq2.3.$$

If $d=2,x_{1}= 5, 2|z,$ then$$x_{1}\frac{\varphi(zd_{1})}{z}\geq 5\times2\frac{\varphi(z)}{z}\geq 5\times2\frac{1}{2}\frac{2}{3}\geq3.3.$$
If $d=2,x_{1}= 5, 2\nshortmid z,z\neq 15,$ then$$x_{1}\frac{\varphi(zd_{1})}{z}\geq 5\frac{\varphi(z)}{z}\geq 5\frac{1}{2}\geq2.5.$$
If $d=2,x_{1}= 5, z=15,$ then$$x_{1}\frac{\varphi(zd_{1})}{z}\geq 5\frac{\varphi(z)}{z}\geq 5\frac{2}{3}\frac{4}{5}\geq2.6.$$
If $d=1,x_{1}=5,z\neq 6$ then $$x_{1}\frac{\varphi(zd_{1})}{z}=5\frac{\varphi(z)}{z}\geq 5\frac{1}{2}\geq2.5.$$
In the above situation, we have $$x_{1}\frac{\varphi(zd_{1})}{z}\geq 2.$$

i)$D_{2}=\phi,$ By Lemma 3.1, for any given prime $q\mid m,$ there are at most 3 primes $p$ with $l_{p}\mid2m, l_{p}\nmid m, q\mid l_{p}.$ It follow that $\sharp\{p: l_{p}\in D_{1}\}\leq9.$
So, $$\prod_{l_{p}\in D_{1}}(1+\frac{1}{p-1})\leq \prod_{4\leq i\leq 13}^{p\neq17}(1+\frac{1}{p-1})<1.72979,$$ where $p_{i}$ is the $i$-th prime.
Thus $2<1.72979,$ a contradiction.

ii)Only two of prime factors of $m$ are less than 173. Then $D_{1}\subseteq \{q,q^{2},q\gamma,q\gamma^{2},q\gamma p\},$ and $D_{2}\subseteq \{\gamma,\gamma^{2},\}.$ By lemma 3.1, for prime $\delta\nmid q\gamma^{3},$ there are at most 3 primes with $l_{p}\mid2m, l_{p}\nmid m, q\mid l_{p};$ It follow that $$\sharp\{p:l_{p}\in D_{1}\}\leq 6.$$

So $$\prod_{l_{p}\in D_{1}}(1+\frac{1}{p-1})\leq \prod_{4\leq i\leq 10}^{p\neq17}(1+\frac{1}{p-1})<1.65.$$

By lemma 3.5, we have $$\prod_{l_{p}\in D_{2}}(1+\frac{1}{p-1})=\prod_{d\in D_{2}}\exp(\log\prod_{l_{p}=d}(1+\frac{1}{p-1})\leq \prod_{d\in D_{2}}\exp0.0355\leq\exp2\times0.03834<1.08.$$
Thus $2<1.65\times1.08=1.782,$ a contradiction.

iii)Only one of prime factor of $m$ are less than 173. Then $D_{1}\subseteq \{q,q^{2},q\gamma,q\gamma^{2},q\gamma p\},$ and $D_{2}\subseteq \{\gamma,\gamma^{2},q\gamma,q\gamma^{2},\}.$ By lemma 4.1, for prime $\delta\nmid q\gamma^{3},$ there are at most 3 primes with $l_{p}\mid2m, l_{p}\nmid m, q\mid l_{p}.$  It follow that $$\sharp\{p:l_{p}\in D_{1}\}\leq 3.$$

So $$\prod_{l_{p}\in D_{1}}(1+\frac{1}{p-1})\leq \prod_{4\leq i\leq 6}(1+\frac{1}{p-1})<1.4.$$

By lemma 3.5, We have $$\prod_{l_{p}\in D_{2}}(1+\frac{1}{p-1})=\prod_{d\in D_{2}}\exp(\log\prod_{l_{p}=d}(1+\frac{1}{p-1}))\leq \prod_{d\in D_{2}}\exp0.03834\leq\exp4\times0.03834<1.17.$$
Thus $2<1.4\times1.17=1.638,$ a contradiction.

Case 2: $d_{1}=1,x_{1}=5,z=6$ then $$x_{1}\frac{\varphi(zd_{1})}{z}=5\frac{\varphi(z)}{z}\geq 5\frac{1}{2}\frac{2}{3}\geq1.667.$$
Noting $3\nmid x_{1}-y_{1}=5-1=4.$ So $m\in \{q,q^{2},q\gamma,q\gamma^{2}\}.$ Through a similar discussion, we get a contradiction.

Case 3: $x=3.$ Then $x_{1}=3, y_{1}=1, d_{1}=1, z\leq x+y\leq4.$ By $q\mid m, 2\nmid m,$

If$q=3,z=3,$ then $q=3\mid m,$ then $m=3^{\nu_{3}(m)}m_{3},$ By lemma 4.2(1) formula (3.2), So
$$\varphi(3)3^{\frac{1}{2}\nu_{3}(m)\tau(m)}\mid 3(3^{n}-1).$$
That is $\frac{1}{2}\nu_{3}(m)\tau(m)=1.$  So $m=3,$
But $$\varphi(3\frac{3^{3}+1}{3+1})=3\frac{3^{n}-1}{3+1}$$ has no solution, a contradiction.

So $z\neq3,$ when $q=3.$  It follow that $q\nmid z.$ By Lemma 3.1(1), $$q^{\frac{1}{2}\nu_{q}(m)\tau(m)}|3^{q-1}-1.$$ A simple calculation shows that, there are no odd primes $p<173$ with $$q^{3}|3^{q-1}-1.$$ So $\frac{1}{2}\nu_{q}(m)\tau(m)\leq2.$ It follow that $m\in \{q,q\gamma\}.$

i)If $m=q\gamma.$ Then $\frac{1}{2}\nu_{q}(m)\tau(m)=2,$ By lemma4.1, $$q^{2}|3^{q-1}-1.$$ Since $q<173,$ it follows from a simple calculation that $q=11.$ this implies that $\gamma\geq 173.$ So $D_{1}=\{22,22\gamma\},D_{2}=\{2\gamma\}.$
Since $3^{11}+1=4\times67\times661,$ it follow that $\{p:l_{p}=22\}=\{67,661\}.$ By Carmichael primitive divisor theorem, $3^{11\gamma}+1$ has at least one primitive prime divisor $p^{\prime}\equiv1(11\gamma).$ By the definition of $l_{p^{\prime}},$ we have $l_{p^{\prime}}=22\gamma.$
Hence $\sharp\{p:l_{p}\in D_{1}\}\geq 3.$

Let $p_{1},p_{2},p_{3}$ be three distinct primes with $l_{p_{i}}\in D_{1}.$ Since
$$p_{i}\mid \frac{3^{l_{p_{i}}}+1}{3+1}, \frac{3^{l_{p_{i}}}+1}{3+1}\mid\frac{3^{m}+1}{3+1},$$ it follow that $$p_{1}p_{2}p_{3}|z\frac{3^{m}+1}{3+1}.$$

So $$(p_{1}-1)(p_{2}-1)(p_{3}-1)|\varphi(z\frac{3^{m}+1}{3+1}).$$

Since $11|l_{p_{i}}, l_{p_{i}}|p_{i}-1,$ It follow that $$11^{3}|\varphi(z\frac{3^{n}-1}{3+1})=z\frac{3^{n}-1}{3+1}.$$
By $z\leq4, 11^{3}|3^{n}-1.$
Noting $3^{5}-1=2\times11^{2},$ we have $11^{2}|3^{(n,5)}-1.$It follows that gcd$(n,5)=5.$ Let $n=5n_{1}.$ Then
$$3^{n}-1=3^{5n_{1}}-1=(2\times11^{2}+1)^{n_{1}}-1=C^{n_{1}}_{1}2\times11^{2}+C^{n_{1}}_{2}(2\times11^{2})^{2}+\cdots+C^{n_{1}}_{n_{1}}(2\times11^{2})^{n_{1}}.$$

By $11^{3}|3^{n}-1,$ we have $11|n_{1},$ then $11|n.$ Since $11=q|m,$ it contradicts gcd$(m,n)=1.$

ii) $m=q.$ Then $D_{1}=\{2q\}, D_{2}=\phi.$ By lemma 4.3 and 5.1, $\sharp\{p: l_{p}=2q\}\leq3.$ Noting that $x_{1}=3, d=1,$ we have
\begin{eqnarray*}
3\frac{\varphi(z)}{z} &<&\prod_{l_{p}\in D_{1}}(1+\frac{1}{p-1})\\
&\leq& (1+\frac{1}{2q})(1+\frac{1}{4q})(1+\frac{1}{6q})\\
&\leq& (1+\frac{1}{6})(1+\frac{1}{12})(1+\frac{1}{18})\\
&<& 1.34\\
\end{eqnarray*}
But for$1\leq z\leq 4,$ $$3\frac{\varphi(z)}{z}\geq 1.5,$$ a contradiction.

This completes the proof of Lemma 3.7. \hfill$\Box$\\

\section{Proof of Theorem 1.4}
$ $

In this section, We always assume that $(x,y,z,m,n)$ is a non trivial solution of the equation (1.6) $$\varphi(z\frac{x^{m}+y^{m}}{x+y})=z\frac{x^{n}-y^{n}}{x+y}, 2\nmid m, 2|n,$$
in positive integers $x,y,z,m,n$ with $1\leq z\leq x+y, z\neq2,$ gcd$(m,n)=1, m>n\geq 1, x>y\geq 1.$

It follow from Theorem 1.3(3) that Theorem 1.4 is true when $\nu_{p}(x)\neq\nu_{p}(y).$ So we always assume that $\nu_{2}(x)=\nu_{2}(y),$ then $x_{1}, y_{1}$ are both odd.

\begin{lemma}\label{lemma1}
$\log x<1.38+\sum_{p|z\frac{x^{m}+y^{m}}{x+y},p\geq 7}\frac{1}{p}.$
\end{lemma}

Proof.By formula (3.1), we have
\begin{eqnarray*}
\log x &<& \sum_{p|z\frac{x^{m}+y^{m}}{x+y}}\log(1+\frac{1}{p-1})\\
&\leq &\log\frac{15}{4}+\sum_{p|z\frac{x^{m}+y^{m}}{x+y},p\geq 7}\log(1+\frac{1}{p-1})\\
&\leq &\log\frac{15}{4}+\sum_{p|z\frac{x^{m}+y^{m}}{x+y},p\geq 7}\frac{1}{p-1}\\
&\leq &\log\frac{15}{4}+\sum_{p\geq 7}\frac{1}{p(p-1)}+\sum_{p|z\frac{x^{m}+y^{m}}{x+y},p\geq 7}\frac{1}{p}\\
&\leq &1.38+\sum_{p|z\frac{x^{m}+y^{m}}{x+y},p\geq 7}\frac{1}{p}\\
\end{eqnarray*}

\begin{lemma}\label{lemma1}
(1)If $p(m)\geq 79,$ then $$\log(\sum_{r\mid m, r>1}\frac{\log\log r}{\varphi(r)}+1)<\frac{\log2p(m)}{\log 2}\frac{\log\log p(m)}{p(m)-1}\frac{p(m)}{p(m)-\log\log p(m)}=f(p(m)).$$
(2)The function $f(x)$ is a monotonically decreasing function.
\end{lemma}

Proof. We follow the proof of [8]. (1) By the lemma 4.8 of [8], If $a, b\geq 78,$ then $\log\log(ab)\leq(\log\log(a))(\log\log(b)),$
\begin{eqnarray*}
\sum_{r|m,r>1}\frac{\log\log r}{\varphi(r)} &<& \prod_{q|m}(1+\frac{\log\log q}{\varphi(q)}+\frac{\log\log q^{2}}{\varphi(q^{2})}+\cdots)-1\\
&\leq& \prod_{q|m}(1+\frac{\log\log q}{\varphi(q)}+\frac{(\log\log q)^{2}}{\varphi(q^{2})}+\cdots)-1\\
&<& \prod_{q|m}(1+\frac{\log\log q}{q-1}\frac{1}{1-\frac{\log\log q}{q}})-1\\
&<& (1+\frac{\log\log p(m)}{p(m)-1}\frac{p(m)}{p(m)-\log\log p(m)})^{\omega(m)}-1\\
\end{eqnarray*}
By lemma3.2, it follow that
\begin{eqnarray*}
\log(1+\sum_{r|m,r>1}\frac{\log\log r}{\varphi(r)}) &<&\omega(m)\log(1+\frac{\log\log p(m)}{p(m)-1}\frac{p(m)}{p(m)-\log\log p(m)})\\
&\leq& \omega(m)\frac{\log\log p(m)}{p(m)-1}\frac{p(m)}{p(m)-\log\log p(m)}\\
&<& \frac{\log2p(m)}{\log 2}\frac{\log\log p(m)}{p(m)-1}\frac{p(m)}{p(m)-\log\log p(m)}\\
&<& \frac{1}{\log 2}\frac{\log2p(m)}{\sqrt{p(m)}}\frac{\log\log p(m)}{\sqrt{p(m)}}\frac{p(m)}{p(m)-1}\frac{p(m)}{p(m)-\log\log p(m)}\\
\end{eqnarray*}

(2) From (1), it is easy to prove.

This completes the proof of Lemma4.2. \hfill$\Box$\\

We divide into three subsections: $ x>73, p(m)\leq x, x>73, p(m)> x$ and $x<73.$

\subsection{$p(m)\leq x, x>73$}

By lemma 4.1, we have
\begin{eqnarray*}
\log x &\leq &1.38+\sum_{p|z\frac{x^{m}+y^{m}}{x+y},p\geq 7}\frac{1}{p}\\
&\leq &1.38+\sum_{p|z\frac{x^{m}+y^{m}}{x+y},7\leq p\leq x^{6}}\frac{1}{p}+\sum_{p|z\frac{x^{m}+y^{m}}{x+y},p>x^{6}}\frac{1}{p}\\
\end{eqnarray*}

By [31], for $t>286,$ $$\sum_{p\leq t}\frac{1}{p}<\log\log t+0.2615+\frac{1}{2\log^{2} t}<\log\log t+0.2772.$$
It follow that
$$\sum_{p|z\frac{x^{m}+y^{m}}{x+y},7\leq p\leq x^{6}}\frac{1}{p}<\log\log x^{6}+0.2772-\frac{1}{2}-\frac{1}{3}-\frac{1}{5}<\log\log x+1.0357$$

Hence
$$\log x<\log\log x +1.38+1.0357+\sum_{p|z\frac{x^{m}+y^{m}}{x+y},p>x^{6}}\frac{1}{p}<\log\log x +2.4157+\sum_{p|z\frac{x^{m}+y^{m}}{x+y},p> x^{6}}\frac{1}{p}.$$

It is clear that, if $p>x^{6},$ then $p\nmid xyz, l_{p}=\min\{l: p|x_{1}^{l}-y_{1}^{l}\}\geq7.$
Hence
$$\sum_{p|z\frac{x^{m}+y^{m}}{x+y},p>x^{6}}\frac{1}{p}=\sum_{p|x_{1}^{m}+y_{1}^{m},p>x^{6}}\frac{1}{p}<\sum_{d|2m,d\nmid m, d\geq7}T_{d},$$
Where $T_{d}=\sum_{l_{p}=d, p>x^{6}}\frac{1}{p}.$

Let $P_{d}=\{p: l_{p}=d,p>x^{6}\},$ then $$x^{6|P_{d}|}<\prod_{p\in P_{d}}p\leq x_{1}^{d}-y_{1}^{d}<x^{d}.$$ It follow that $|P_{d}|<\frac{d}{6}.$
If $d\leq x^{3},$ then $$T_{d}\leq\frac{|P_{d}|}{x^{6}}<\frac{d}{6x^{6}}\leq\frac{1}{6x^{3}}.$$
Thus $$\sum_{d|2m, d\nmid m, 7\leq d\leq x^{3}}T_{d}<x^{3}\frac{1}{6x^{3}}=0.1667.$$

If $d>x^{3}>73^{3}.$ By[25], Let $N\geq 3,$ then $\frac{N}{\varphi(N)}\leq1.79\log\log N+\frac{2.5}{\log\log N};$
By lemma 3.3, we have
\begin{eqnarray*}
T_{d}\leq S_{d} &<&\frac{1.084}{d}+\frac{1}{d\log(d+1)}+\frac{2\log\log d}{\varphi(d)}+\frac{2\log\log x}{\varphi(d)\log d}\\
&\leq&\frac{1.084}{d}+\frac{1}{d\log(d+1)}+\frac{3.58(\log\log d)^{2}}{d}\\
&+&\frac{5}{d}+\frac{3.58\log\log d\log\log x}{d\log d}+\frac{5\log\log x}{d(\log d)\log\log d}\\
&\leq&\frac{1.084}{x^{3}}+\frac{1}{x^{3}\log(x^{3}+1)}+\frac{3.58(\log\log x^{3})^{2}}{x^{3}}\\
&+&\frac{5}{x^{3}}+\frac{3.58\log\log x^{3}\log\log x}{x^{3}\log x^{3}}+\frac{5\log\log x}{x^{3}(\log x^{3})\log\log x^{3}}\\
\end{eqnarray*}

By Lemma 3.2 and $p(m)\leq x,$ we have $\tau(m)<2x,$ So $\sharp\{d: d|2m, d\nmid m, d>x^{3}\}\leq\tau(m)<2x.$ Hence, we have
\begin{eqnarray*}
\sum_{d|2m, d\nmid m, d>x^{3}}T_{d} &< &\frac{2.168}{x^{2}}+\frac{2}{x^{2}\log(x^{3}+1)}+\frac{7.16(\log\log x^{3})^{2}}{x^{2}}\\
&+&\frac{10}{x^{2}}+\frac{7.16\log\log x^{3}\log\log x}{x^{2}\log x^{3}}+\frac{10\log\log x}{x^{2}(\log x^{3})\log\log x^{3}}\\
&<&0.04\\
\end{eqnarray*}

Therefore
\begin{eqnarray*}
\log x &< &\log\log x +2.4157+\sum_{p|z\frac{x^{m}+y^{m}}{x+y},p> x^{6}}\frac{1}{p}\\
&<&\log\log x +2.4157+\sum_{d|2m, d\nmid m, d\geq 7}T_{d}\\
&=&\log\log x +2.4157+\sum_{d|2m, d\nmid m, 7\leq d\leq x^{3}}T_{d}+\sum_{d|2m,d\nmid m, d> x^{3}}T_{d}\\
&<&\log\log x +2.4157+0.1667+0.04\\
&<&\log\log x +2.63\\
\end{eqnarray*}

Since $x>73,$ it follows that $\log x-\log\log x>\log73-\log\log73>2.83.$ a contradiction.

\subsection{$p(m)>x>73$}

Since $p(m)p(m)>x>73,$ it follow that $p(m)\geq79.$ By lemma 4.1, We have
\begin{eqnarray*}
\log x &\leq &1.38+\sum_{p|z\frac{x^{m}+y^{m}}{x+y},p\geq 7}\frac{1}{p}\\
&\leq &1.38+\sum_{p|z\frac{x^{m}+y^{m}}{x+y},7\leq p\leq x^{2}}\frac{1}{p}+\sum_{p|z\frac{x^{m}+y^{m}}{x+y},p>x^{2}}\frac{1}{p}\\
\end{eqnarray*}

By [31], for $t>286,$ $$\sum_{p\leq t}\frac{1}{p}<\log\log t+0.2615+\frac{1}{2\log^{2} t}.$$ It follow from $x>73$ that
$$\sum_{7\leq p\leq x^{2}}\frac{1}{p}<\log\log x^{2}+0.2615+\frac{1}{2\log^{2} 5329}-\frac{1}{2}-\frac{1}{3}-\frac{1}{5}<\log\log x-0.07.$$

It is clear that, if $p>x^{2},$ then $p\nmid xyz.$  $l_{p}=\min\{l: p|x_{1}^{l}-y_{1}^{l}\}>2.$  Hence
\begin{eqnarray*}
\log x &<& 1.38+\log\log x-0.07++\sum_{p|z\frac{x^{m}+y^{m}}{x+y},p>x^{2}}\frac{1}{p}\\
&\leq&\log\log x+1.31+\sum_{p|x_{1}^{m}+y_{1}^{m},p>x^{2}}\frac{1}{p}\\
&\leq&\log\log x+1.31+\sum_{d|2m, d\nmid m, d>2}S_{d}\\
\end{eqnarray*}
where $S_{d}=\sum_{l_{p}=d}\frac{1}{p}.$

For $d|2m, d\nmid m, d>2,$ we have $d\geq 2p(m), p(m)>x>73,$ that is $d=2r\geq 158, 2\nmid r.$
By lemma4.5,
\begin{eqnarray*}
S_{d} &<& \frac{1.084}{d}+\frac{1}{d\log(d+1)}+\frac{2\log\log d}{\varphi(d)}+\frac{2\log\log x}{\varphi(d)\log d}\\
&<& \frac{1.084}{d}+\frac{1}{d\log159}+\frac{2\log\log d}{\varphi(d)}+\frac{2\log\log x}{\varphi(d)\log2x}\\
&<& \frac{1.084}{\varphi(d)}+\frac{1}{\varphi(d)\log159}+\frac{2\log\log d}{\varphi(d)}+\frac{2\log\log120}{\varphi(73)\log146}\\
&<& \frac{1.8659}{\varphi(d)}+\frac{2\log\log d}{\varphi(d)}\\
&<& \frac{3.16\log\log d}{\varphi(d)}\\
&<& \frac{4\log\log r}{\varphi(r)}
\end{eqnarray*}

Since $p(m)>73,$ that is $p(m)\geq79,$ it follow from Lemma 4.2 that
\begin{eqnarray*}
\log(1+\sum_{r|m,r>1}\frac{\log\log r}{\varphi(r)}) &<&\frac{\log2p(m)}{\log 2}\frac{\log\log p(m)}{p(m)-1}\frac{p(m)}{p(m)-\log\log p(m)}\\
&\leq& \frac{\log158}{\log2}\frac{\log\log79}{78}\frac{79}{79-\log\log79}<0.15\\
\end{eqnarray*}
$$\sum_{d|2m, d\nmid m, d>2}S_{d}<4\sum_{r|m,r>1}\frac{\log\log r}{\varphi(r)}<4\times(\exp(0.15)-1)<0.68.$$

So $$\log x<\log\log x+1.31+\sum_{d|2m, d\nmid m, d>2}S_{d}<\log\log x+1.31+0.68<\log\log x+2.$$
Since $x>73,$ it follow that $\log x-\log\log x>2.8,$ a contradiction.

\subsection{$x\leq73$}

By lemma 3.7, $p(m)\geq173.$ It follow that $p(m)>x.$ For $d\mid 2m, d\nmid m, d>2,$ we have $d\geq 2p(m)\geq346.$ In view of Lemma 3.4, 3.5 We have
\begin{eqnarray*}
\log(x_{1}\frac{\varphi(zd_{1})}{z})&<&\log(\prod_{l_{p}\mid2m,l_{p}\nmid m,l_{p}>2}(\frac{p}{p-1}))\\
&\leq& \sum_{d\mid2m,d\nmid m,d>2}\log(\prod_{l_{p}=d}(\frac{p}{p-1}))\\
&\leq& 3.7341\sum_{d\mid2m,d\nmid m,d>2}\frac{\log\log d}{\varphi(d)}\\
&\leq& 4.4903\sum_{r\mid m,r>1}\frac{\log\log r}{\varphi(r)}\\
\end{eqnarray*}
where $d=2r.$

By lemma 3.2 and $p(m)>x,$ $\omega(m)<\frac{\log(2p(m))}{\log2}.$ In view of lemma 4.2, we have
\begin{eqnarray*}
\log(1+\sum_{r|m,r\geq 1}\frac{\log\log r}{\varphi(r)})&<&\frac{\log2p(m)}{\log 2}\frac{\log\log p(m)}{p(m)-1}\frac{p(m)}{p(m)-\log\log p(m)}\\
&\leq& \frac{\log346}{\log2}\frac{\log\log 173}{172}\frac{173}{173-\log\log 173}<0.082\\
\end{eqnarray*}

It follow that
$$\log(x_{1}\frac{\varphi(zd_{1})}{z})<4.4903\sum_{r\mid m,r>1}\frac{\log\log r}{\varphi(r)}<4.4903\times(\exp 0.082-1)<0.384.$$
Hence $$x_{1}\frac{\varphi(zd_{1})}{z}<\exp0.3833<1.47.$$

Since $2\nmid x_{1}y_{1}, x_{1}>y_{1}\geq1,$ it follow that $x_{1}\geq3, 1\leq z\leq120.$ If $d\geq 3,$ then
$$x_{1}\frac{\varphi(zd_{1})}{z}\geq 3\varphi(d_{1})\frac{\varphi(z)}{z}\geq 3\times2\frac{1}{2}\frac{2}{3}\frac{4}{5}\geq1.6.$$

If $d_{1}=2,1,x_{1}\geq7,$ then $$x_{1}\frac{\varphi(zd_{1})}{z}\geq 7\frac{\varphi(z)}{z}\geq 7\frac{1}{2}\frac{2}{3}\frac{4}{5}\geq1.87.$$

If $d_{1}=2,1,x_{1}=5,$ then $$x_{1}\frac{\varphi(zd_{1})}{z}\geq 5\frac{\varphi(z)}{z}\geq 7\frac{1}{2}\frac{2}{3}\geq1.66.$$

If $d_{1}=2,x_{1}= 3, z\neq6,$ then$$x_{1}\frac{\varphi(zd_{1})}{z}\geq 3\frac{\varphi(z)}{z}\geq 3\frac{1}{2}\geq1.5.$$
If $d_{1}=2,x_{1}= 3, z=6,$ then$$x_{1}\frac{\varphi(zd_{1})}{z}= 3\frac{\varphi(12)}{12}=2.$$
If $d_{1}=1,x_{1}= 3,$ then$$x_{1}\frac{\varphi(zd_{1})}{z}\geq 3\frac{\varphi(z)}{z}\geq 3\frac{1}{2}\geq1.5.$$

In any way, we have $$x_{1}\frac{\varphi(zd_{1})}{z}\geq1.5.$$ A contradiction.

This completes the proof of Theorem 1.4. \hfill$\Box$\\

\section{Acknowledgments}

Sincere thanks to Professor Pingzhi Yuan for his careful guidance. This work was supported by the National Nature Science Foundatin of China, No.11671153, No.11971180.

 \end{document}